\documentclass[12pt,a4paper,reqno]{amsart}
\usepackage{amssymb}
\usepackage{dcpic,pictex}
\usepackage{multirow}
\usepackage{graphicx}
\usepackage{slashbox}
\usepackage{cite}

\usepackage{hyperref}
\usepackage[main=english,french]{babel}
\newtheorem{definition}{Definition}[section]

\newtheorem{theorem}{Theorem}[section]
\newtheorem{remark}{Remark}[section]
\newcommand{\D}{\mathbb{D}}

\newcommand{\I}{\mathbb{I}}
\newcommand{\R}{\mathbb{R}}

\newcommand{\Com}{\mathbb{C}}
\newcommand{\N}{\mathbb{N}}

\begin{document}


\title[Special Functions of Fractional Calculus]{Special Functions of Fractional Calculus in Form of Convolution Series and their Applications}

\author{Yuri Luchko}
\curraddr{Beuth Technical University of Applied Sciences Berlin,  
     Department of  Mathematics, Physics, and Chemistry,  
     Luxemburger Str. 10,  
     13353 Berlin,
Germany}
\email{luchko@beuth-hochschule.de}

\subjclass[2010]{26A33; 26B30; 44A10; 45E10}
\dedicatory{}
\keywords{Sonine kernel; Sonine condition; general fractional derivative; general fractional integral; convolution series; fundamental theorems of Fractional Calculus; fractional differential equations}

\begin{abstract}
In this paper, we first discuss the convolution series that are generated by the Sonine kernels from a class of functions continuous on the real positive semi-axis that have an integrable singularity of power function type at the point zero. These convolution series are closely related to the general fractional integrals and derivatives with the Sonine kernels and represent a new class of the special functions of Fractional Calculus. The Mittag-Leffler functions as solutions to the fractional differential equations with the fractional derivatives of both Riemann-Liouville and Caputo types are particular cases of the convolution series generated by the Sonine kernel $\kappa(t) = t^{\alpha -1}/\Gamma(\alpha),\ 0<\alpha <1$. The main result of the paper is derivation of analytic solutions to the single- and multi-term fractional differential equations with the general fractional  derivatives of the Riemann-Liouville type that were not yet considered in the Fractional Calculus publications. 
\end{abstract}

\maketitle


\section{Introduction}
\label{sec:1}

Special functions of mathematical physics are usually defined either in form of power series, or as solutions to some differential equations, or via integral representations. Of course, for a given function, these three (and possibly other) forms coincide for all arguments and parameters values for that they do exist. However, the validity domains of different representations can be unequal. Very often the series representations of the special functions hold valid only on some restricted domains. In order to define the corresponding functions for other values of their arguments and parameters, analytical continuation of the series in form of integral representations is usually employed.  

For the special functions of Fractional Calculus (FC), the situation is very similar to the one described above. For instance, one of the most important FC special function - the two-parameters Mittag-Leffler function - is usually defined in form of a power series:
\begin{equation}
\label{ML}
E_{\alpha,\beta}(z) = \sum_{k=0}^{+\infty} \frac{z^k}{\Gamma(\alpha\, k + \beta)},\ \alpha >0,\ \beta,z\in \Com.
\end{equation}
Because the series is convergent for all $z\in \Com$, this definition can be used for all $z\in \Com$ without any analytical continuation. Still, the integral representations of the  Mittag-Leffler function are very important, say, for derivation of its asymptotic behavior (\cite{Dzh}) and for its numerical calculation (\cite{GLL}). For 
$0<\alpha <2$ and $\Re(\beta)>0$,  the following integral representations of the Mittag-Leffler function in terms of the integrals over the Hankel-type contours were presented in \cite{Dzh}: 
$$
E_{\alpha,\beta}(z) ={1 \over 2\pi\alpha i}\int_{\gamma(\epsilon;\delta)}
{e^{\zeta^{1/\alpha}}\zeta^{(1-\beta)/\alpha} \over \zeta -z}\, d\zeta, \ z\in G^{(-)}(\epsilon;\delta),
$$
$$
E_{\alpha,\beta}(z) = \frac{1}{\alpha}z^{(1-\beta)/\alpha}e^{z^{1/\alpha}} + {1\over 2\pi\alpha i}\int_{\gamma(\epsilon;\delta)}
{e^{\zeta^{1/\alpha}}\zeta^{(1-\beta)/\alpha} \over \zeta -z} d\zeta, \ z\in G^{(+)}(\epsilon;\delta),
$$
where   the integration contour $\gamma(\epsilon;\delta) \ (\epsilon >0, \ 0<\delta\le \pi)$  with
non-decreasing $\arg \zeta$ consists of the following parts:
\par \noindent
1) the ray $\arg \zeta =-\delta, \ |\zeta|\ge \epsilon$;
\par \noindent
2) the arc $-\delta \le \arg \zeta \le \delta$ of the circumference $|\zeta|=\epsilon$;
\par \noindent
3) the ray  $\arg \zeta =\delta, \ |\zeta|\ge \epsilon$.
\par \noindent
For $0<\delta< \pi$, the domain $G^{(-)}(\epsilon;\delta)$ is to the left
 of the contour $\gamma(\epsilon;\delta)$ and the domain  $G^{(+)}(\epsilon;\delta)$ is to the right of this contour. If $\delta
=\pi $, the contour $\gamma(\epsilon;\delta)$  consists of the circumference $|\zeta|=\epsilon $ and of the cut
$-\infty <\zeta \le -\epsilon$. In this case, the domain  $G^{(-)}(\epsilon;\delta)$ is the circle $|\zeta|<
\epsilon$ and  $G^{(+)}(\epsilon;\alpha)=\{\zeta: |\arg\zeta| <\pi, \ |\zeta|>\epsilon \}$.

For some parameters values, the Mittag-Leffler function can be also introduced in terms of solutions to the fractional differential equations with the Riemann-Liouville or Caputo fractional derivatives. For instance, for $0<\alpha \le 1$, the equation 
\begin{equation}
\label{eqRL}
(D^\alpha_{0+} y)(t)  = \lambda\, y(t)
\end{equation}
has the general solution (\cite{LucSri})
\begin{equation}
\label{eqRL_sol}
y(t) = C\, t^{\alpha -1} E_{\alpha, \alpha}(\lambda\, t^\alpha),\ C\in \R.
\end{equation}
In the equation \eqref{eqRL}, the Riemann-Liouville fractional derivative $D^\alpha_{0+}$ is defined by 
\begin{equation}
\label{RLD}
(D^\alpha_{0+} \, f)(t) = \frac{d}{dt}(I^{1-\alpha}_{0+} f)(t),\ t>0,
\end{equation}
where $I^\alpha_{0+}$ is the Riemann-Liouville fractional integral of order $\alpha\ (\alpha >0)$:
\begin{equation}
\label{RLI}
\left(I^\alpha_{0+} f \right) (t)=\frac{1}{\Gamma (\alpha )}\int\limits_0^t(t -\tau)^{\alpha -1}f(\tau)\,d\tau,\ t>0.
\end{equation}
The general solution to the equation 
\begin{equation}
\label{eqC}
(\,_*D^\alpha_{0+} y)(t)  = \lambda\, y(t)
\end{equation}
with the Caputo fractional derivative
\begin{equation}
\label{CD}
(\,_*D^\alpha_{0+}\, f)(t)  = (D^\alpha_{0+} \, f)(t) - f(0)\frac{t^{-\alpha}}{\Gamma(1-\alpha)},\ t>0 
\end{equation}
has the form  (\cite{LucGor})
\begin{equation}
\label{eqC_sol}
y(t) = C\,  E_{\alpha, 1}(\lambda\, t^\alpha),\ C\in \R.
\end{equation}
As we see, the solutions to the fractional differential equations \eqref{eqRL} and \eqref{eqC} are expressed in terms of the Mittag-Leffler functions. However, the arguments of these functions are  $\lambda\, t^\alpha$ and not just $\lambda\, t$. Thus, these solutions are represented in form of power series with the fractional and not integer exponents. For more advanced properties and applications of the Mittag-Leffler type functions see \cite{Dzh} and the recent book \cite{GKRM}.  

In \cite{Luc21b}, the single- and multi-terms fractional differential equations with the general fractional derivatives of Caputo type have been studied. By definition, their solutions belong to the class of the FC special functions (as the ones represented in form of solutions to the fractional differential equations). Moreover, in \cite{Luc21b}, another representation of these new FC special functions was derived, namely, in form of the convolution series generated by the Sonine kernels.  

The convolution series are a far reaching generalization of the conventional power series and the power series with the fractional exponents including the Mittag-Leffler functions \eqref{eqRL_sol} and \eqref{eqC_sol}. They represent a new type of the FC special functions worth for investigation. In \cite{Luc21d}, the convolution series were employed for derivation of two different forms of a generalized convolution Taylor formula for representation of a function as a convolution polynomial with a remainder in form of a composition of the $n$-fold general fractional integral and the $n$-fold general sequential fractional derivative of the Riemann-Liouville and the Caputo types, respectively.  In \cite{Luc21d}, the generalized Taylor series in form of convolution series were also discussed. In this paper, we employ the convolution series for derivation of analytical solutions to the single- and multi-terms fractional differential equations with the general fractional derivatives in the Riemann-Liouville  sense.  This type of the fractional differential equations was not yet considered in the FC literature.

The rest of the paper is organized as follows. In the 2nd section, we introduce the general fractional derivatives of the Riemann-Liouville and Caputo types with the Sonine kernels from a special class of kernels and discuss some of their properties needed for the further discussions. In the 3rd section, we first provide some results  regarding the convolution series generated by the Sonine kernels. Then the convolution series are applied for derivation of analytical solutions to the single- and multi-terms fractional differential equations with the general fractional derivatives in the Riemann-Liouville  sense. For a treatment of the single- and multi-terms fractional differential equations with the general fractional derivatives in the Caputo  sense, we refer the interested readers to \cite{Luc21b}.

\section{General Fractional Integrals and Derivatives}

The general fractional derivatives (GFDs) with the kernel $k$ in the Riemann-Liouville and in the Caputo sense, respectively, are defined as follows (\cite{Koch11,LucYam16,Koch19_1,LucYam20,Luc21a,Luc21c}): 
\begin{equation}
\label{FDR-L}
(\D_{(k)}\, f) (t) = \frac{d}{dt}(k\, *\, f)(t)= \frac{d}{dt}\int_0^t k(t-\tau)f(\tau)\, d\tau,
\end{equation}
\begin{equation}
\label{FDC}
( _*\D_{(k)}\, f) (t) =  (\D_{(k)}\, f) (t) - f(0)k(t),
\end{equation}
where by $*$ the Laplace convolution is denoted:
\begin{equation}
\label{2-2}
(f\, *\, g)(t) = \int_0^{t}\, f(t-\tau)g(\tau)\, d\tau.
\end{equation}

The Riemann-Liouville and the Caputo  fractional derivatives of order $\alpha,\ 0< \alpha < 1$, defined by \eqref{RLD} and \eqref{CD}, respectively, are particular cases of the GFDs \eqref{FDR-L} and \eqref{FDC} with the kernel
\begin{equation}
\label{single}
k(t) = h_{1-\alpha}(t),\ 0 <\alpha <1,\ h_{\beta}(t) := \frac{t^{\beta -1}}{\Gamma(\beta)},\ \beta >0.
\end{equation}

The multi-term fractional derivatives and the fractional derivatives of the
distributed order are also particular cases of the GFDs \eqref{FDR-L} and \eqref{FDC} with the kernels 
\begin{equation}
\label{multi}
k(t) = \sum_{k=1}^n a_k\, h_{1-\alpha_k}(t),
\ \  0 < \alpha_1 <\dots < \alpha_n < 1,\ a_k \in \R,\ k=1,\dots,n, 
\end{equation}
\begin{equation}
\label{distr}
k(t) = \int_0^1 h_{1-\alpha}(t)\, d\rho(\alpha),
\end{equation}
respectively, where $\rho$ is a Borel measure defined on the interval $[0,\, 1]$.

Several useful properties of the Riemann-Liouville fractional integral and the Riemann-Liouville and Caputo fractional derivatives are based on the formula
\begin{equation}
\label{2-9}
(h_{\alpha} \, * \, h_\beta)(t) \, = \, h_{\alpha+\beta}(t),\ \alpha,\beta >0,\ t>0
\end{equation}
that immediately follows from the well-known representation of the Euler Beta-function in terms of the Gamma-function:
$$
B(\alpha,\beta) := \int_0^1 (1-\tau)^{\alpha -1}\, \tau^{\beta -1}\, d\tau \, = \, \frac{\Gamma(\alpha)\Gamma(\beta)}{\Gamma(\alpha+\beta)},\ \alpha,\beta>0.
$$ 
In the formula \eqref{2-9} and in what follows,  the power function $h_{\alpha}$ is defined as in \eqref{single}. 

In our discussions, we employ the integer order convolution powers that for a function $f=f(t),\ t>0$ are defined by the expression
\begin{equation}
\label{I-6}
f^{<n>}(t):= \begin{cases}
1,& n=0,\\
f(t), & n=1,\\
(\underbrace{f* \ldots * f}_{n\ \mbox{times}})(t),& n=2,3,\dots .
\end{cases}
\end{equation}
For the kernel $\kappa(t) = h_\alpha(t)$ of the Riemann-Liouville fractional integral, we apply the formula \eqref{2-9} and arrive at the important representation
\begin{equation}
\label{2-13}
h_{\alpha}^{<n>}(t) = h_{n\alpha}(t),\ n\in \N.
\end{equation}
A well-known particular case of \eqref{2-13} is the formula
\begin{equation}
\label{2-13-1}
\{1\}^n(t) = h_{1}^n(t)= h_{n}(t)=\frac{t^{n-1}}{\Gamma(n)} = \frac{t^{n-1}}{(n-1)!},\ n\in \N,
\end{equation}
where by $\{1\}$ we denoted a function that is identically equal to 1 for $t\ge 1$. 

Now let us write down the formula \eqref{2-9} for $\beta = 1-\alpha,\ 0<\alpha <1$: 
\begin{equation}
\label{3-1}
(h_{\alpha}\, * \, h_{1-\alpha})(t) = h_1(t) = \{1 \},\ 0<\alpha<1,\  t>0.
\end{equation}

In \cite{Abel1,Abel2}, Abel employed the relation \eqref{3-1} to derive an inversion formula for the operator that is nowadays referred to as the Caputo fractional derivative and obtained it in form of the Riemann-Liouville fractional integral (solution to the Abel model for the tautochrone problem). 

By an attempt to extend the Abel solution  method to more general integral equations of convolution type, 
Sonine introduced in \cite{Son} the relation
\begin{equation}
\label{3-2}
(\kappa \, *\, k )(t) = \{1 \},\  t>0
\end{equation}
that  is nowadays referred to as the Sonine condition. The functions that satisfy the Sonine condition are called the Sonine kernels. For a Sonine kernel $\kappa$, the kernel $k$ that satisfies the Sonine condition \eqref{3-2} is called an associated kernel to $\kappa$. Of course, $\kappa$ is then an associated kernel to $k$. In what follows, we denote the set of the Sonine kernels by $\mathcal{S}$.   

In \cite{Son}, Sonine introduced a class of the Sonine kernels in the form
\begin{equation}
\label{3-3}
\kappa(t) = h_{\alpha}(t) \cdot \, \kappa_1(t),\ \kappa_1(t)=\sum_{k=0}^{+\infty}\, a_k t^k, \ a_0 \not = 0,\ 0<\alpha <1,
\end{equation}
\begin{equation}
\label{3-4}
k(t) = h_{1-\alpha}(t) \cdot k_1(t),\ k_1(t)=\sum_{k=0}^{+\infty}\, b_k t^k, 
\end{equation}
where $\kappa_1=\kappa_1(t)$ and $k_1=k_1(t)$ are analytical functions and the coefficients $a_k,\ b_k,\ k\in \N_0$  satisfy the following triangular system of linear equations:
\begin{equation}
\label{3-5}
a_0b_0 = 1,\ \sum_{k=0}^n\Gamma(k+1-\alpha)\Gamma(\alpha+n-k)a_{n-k}b_k = 0,\ n\ge 1.
\end{equation}
An important example of the kernels from $\mathcal{S}$ in form \eqref{3-3}, \eqref{3-4} was derived in \cite{Son}  in terms of the Bessel function $J_{\nu}$ and the modified Bessel function $I_{\nu}$: 
\begin{equation}
\label{Bess}
\kappa(t) = (\sqrt{t})^{\alpha-1}J_{\alpha-1}(2\sqrt{t}),\ 
k(t) = (\sqrt{t})^{-\alpha}I_{-\alpha}(2\sqrt{t}),\ 0<\alpha <1,
\end{equation}
where 
$$
J_\nu (t) = \sum_{k=0}^{+\infty} \frac{(-1)^k(t/2)^{2k+\nu}}{k!\Gamma(k+\nu+1)},\ 
I_\nu (t) = \sum_{k=0}^{+\infty} \frac{(t/2)^{2k+\nu}}{k!\Gamma(k+\nu+1)}.
$$

For other examples of the Sonine kernels we refer the readers to \cite{Luc21c,Koch11,Luc21a,Sam}.

In this paper, we deal with the general fractional integrals (GFIs) with the kernels $\kappa \in \mathcal{S}$ 
defined by the formula
\begin{equation}
\label{GFI}
(\I_{(\kappa)}\, f)(t) := (\kappa\, *\, f)(t) = \int_0^t \kappa(t-\tau)f(\tau)\, d\tau,\ t>0
\end{equation}
and with the GFDs with the associated Sonine kernels $k$ in the Riemann-Liouville and Caputo senses defined by \eqref{FDR-L} and \eqref{FDC}, respectively.

In our discussions, we restrict ourselves to a class of the Sonine kernels from the space $C_{-1,0}(0,+\infty)$ that is an important particular case of the following two-parameters family of spaces (\cite{Luc21b,Luc21a,Luc21c}):
\begin{equation}
\label{subspace}
 C_{\alpha,\beta}(0,+\infty) \, = \, \{f:\ f(t) = t^{p}f_1(t),\ t>0,\ \alpha < p < \beta,\ f_1\in C[0,+\infty)\}.
\end{equation}
By $C_{-1}(0,+\infty)$ we mean the space $C_{-1,+\infty}(0,+\infty)$.

The set of such Sonine kernels will be denoted by $\mathcal{L}_{1}$ (\cite{Luc21c}):
\begin{equation}
\label{Son}
(\kappa,\, k \in \mathcal{L}_{1} ) \ \Leftrightarrow \ (\kappa,\, k \in C_{-1,0}(0,+\infty))\wedge ((\kappa\, *\, k)(t) \, = \, \{1\}).
\end{equation}

In the rest of this section, we present some important results for the GFIs and the GFDs  with the Sonine kernels from $\mathcal{L}_{1}$ on the space $C_{-1}(0,+\infty)$ and its sub-spaces. 

The basic properties of the GFI \eqref{GFI} on the space $C_{-1}(0,+\infty)$ easily follow from the known properties of the Laplace convolution:
\begin{equation}
\label{GFI-map}
\I_{(\kappa)}:\, C_{-1}(0,+\infty)\, \rightarrow C_{-1}(0,+\infty),
\end{equation}
\begin{equation}
\label{GFI-com}
\I_{(\kappa_1)}\, \I_{(\kappa_2)} = \I_{(\kappa_2)}\, \I_{(\kappa_1)},\ \kappa_1,\, \kappa_2 \in \mathcal{L}_{1},
\end{equation}
\begin{equation}
\label{GFI-index}
\I_{(\kappa_1)}\, \I_{(\kappa_2)} = \I_{(\kappa_1*\kappa_2)},\  \kappa_1,\, \kappa_2 \in \mathcal{L}_{1}.
\end{equation}

For the functions $f\in C_{-1}^1(0,+\infty):=\{f:\ f^\prime\in C_{-1}(0,+\infty)\}$, the GFDs of the Riemann-Liouville type can be represented as follows (\cite{Luc21a}):
\begin{equation}
\label{GFDL-1}
(\D_{(k)}\, f) (t) = (k\, * \, f^\prime)(t) + f(0)k(t),\ t>0.
\end{equation}
Thus, for $f\in C_{-1}^1(0,+\infty)$, the GFD \eqref{FDC} of the Caputo type takes the form
\begin{equation}
\label{GFDC_1}
( _*\D_{(k)}\, f) (t) =  (k\, * \, f^\prime)(t),\ t>0. 
\end{equation}

It is worth mentioning that in the FC publications, the Caputo fractional derivative \eqref{CD} is often defined as in the formula \eqref{GFDC_1}:
\begin{equation}
\label{CD-1}
(\,_*D^\alpha_{0+}\, f)(t)  = (h_{1-\alpha}\, * \, f^\prime)(t) = (I^{1-\alpha}_{0+} f^\prime)(t),\ t>0. 
\end{equation}

Now, following \cite{Luc21a,Luc21d}, we define the $n$-fold GFI and the $n$-fold sequential GFDs in the Riemann-Liouville and Caputo sense. 

\begin{definition}[\cite{Luc21a}]
\label{d1}
Let $\kappa \in \mathcal{L}_{1}$. The $n$-fold GFI ($n \in \N$) is  a composition of $n$ GFIs with the kernel $\kappa$:
\begin{equation}
\label{GFIn}
(\I_{(\kappa)}^{<n>}\, f)(t) := (\underbrace{\I_{(\kappa)} \ldots \I_{(\kappa)}}_{n\ \mbox{times}}\, f)(t),\  t>0.
\end{equation}
\end{definition}

It is worth mentioning that  the index law \eqref{GFI-index} leads to a representation of the $n$-fold GFI \eqref{GFIn} in form of the GFI with the kernel $\kappa^{<n>}$:
\begin{equation}
\label{GFIn-1}
(\I_{(\kappa)}^{<n>}\, f)(t) = (\kappa^{<n>}\, *\, f)(t) = (\I_{(\kappa)^{<n>}}\, f)(t),\  t>0.
\end{equation}

The kernel $\kappa^{<n>},\ n\in \N$ belongs to the space $C_{-1}(0,+\infty)$, but it is not always a Sonine kernel. 

\begin{definition}[\cite{Luc21d}]
\label{d2}
Let $\kappa \in \mathcal{L}_{1}$ and $k$ be its associated Sonine kernel. The $n$-fold sequential GFDs in the Riemann-Liouville and in the Caputo sense, respectively, are defined as follows:
\begin{equation}
\label{GFDLn}
(\D_{(k)}^{<n>}\, f)(t) := (\underbrace{\D_{(k)} \ldots \D_{(k)}}_{n\ \mbox{times}}\, f)(t),\  t>0,
\end{equation}
\begin{equation}
\label{GFDLn-C}
(\,_*\D_{(k)}^{<n>}\, f)(t) := (\underbrace{\,_*\D_{(k)} \ldots \,_*\D_{(k)}}_{n\ \mbox{times}}\, f)(t),\  t>0.
\end{equation}
\end{definition}

It is worth mentioning that in \cite{Luc21b,Luc21a}, the $n$-fold  GFDs ($n \in \N$) were defined in a different form:
\begin{equation}
\label{GFDLn-1}
(\D_{(k)}^n\, f)(t) := \frac{d^n}{dt^n} ( k^{<n>} * f)(t),\ t>0,
\end{equation}
\begin{equation}
\label{GFDLn-1_C}
(\,_*\D_{(k)}^n\, f)(t) := ( k^{<n>} * f^{(n)})(t),\ t>0.
\end{equation}

The $n$-fold sequential GFDs \eqref{GFDLn} and \eqref{GFDLn-C} are a far reaching generalization of the Riemann-Liouville and the Caputo sequential fractional derivatives to the case of the Sonine kernels from $\mathcal{L}_{1}$.

Some important connections between the $n$-fold GFI \eqref{GFIn} and the $n$-fold sequential GFDs \eqref{GFDLn} and \eqref{GFDLn-C} in the Riemann-Liouville and Caputo senses are provided in the so-called first and second fundamental theorems of FC (\cite{Luc20}) formulated below.  

\begin{theorem}[\cite{Luc21d}]
\label{t3-n}
Let $\kappa \in \mathcal{L}_{1}$ and $k$ be its associated Sonine kernel. 

Then,  the $n$-fold sequential GFD \eqref{GFDLn} in the Riemann-Liouville sense  is a left inverse operator to the $n$-fold  GFI \eqref{GFIn} on the space $C_{-1}(0,+\infty)$: 
\begin{equation}
\label{FTLn}
(\D_{(k)}^{<n>}\, \I_{(\kappa)}^{<n>}\, f) (t) = f(t),\ f\in C_{-1}(0,+\infty),\ t>0,
\end{equation}
and the $n$-fold sequential GFD \eqref{GFDLn-C} in the Caputo sense  is a left inverse operator to the $n$-fold  GFI \eqref{GFIn} on the space $C_{-1,(k)}^n(0,+\infty)$: 
\begin{equation}
\label{FTLn-C}
(\,_*\D_{(k)}^{<n>}\, \I_{(\kappa)}^{<n>}\, f) (t) = f(t),\ f\in C_{-1,(k)}^n(0,+\infty),\ t>0,
\end{equation}
where $C_{-1,(k)}^n(0,+\infty) := \{f:\ f(t)=(\I_{(k)}^{<n>}\, \phi)(t),\ \phi\in C_{-1}(0,+\infty)\}$.
\end{theorem}

\begin{theorem}[\cite{Luc21d}]
\label{tgcTfn}
Let $\kappa \in \mathcal{L}_{1}$ and $k$ be its associated Sonine kernel. 

For a function $f\in C_{-1,(k)}^{(n)}(0,+\infty) = \{ f:\ (\D_{(k)}^{<j>}\, f)\in C_{-1}(0,+\infty),\ j=0,1,\dots,n\}$, the formula
\begin{equation}
\label{sFTLn}
(\I_{(\kappa)}^{<n>}\, \D_{(k)}^{<n>}\, f) (t) = f(t) - \sum_{j=0}^{n-1}\left( k\, * \, \D_{(k)}^{<j>}\, f\right)(0)\kappa^{<j+1>}(t) = 
\end{equation}
$$
f(t) - \sum_{j=0}^{n-1}\left( \I_{(k)}\, \D_{(k)}^{<j>}\, f\right)(0)\kappa^{<j+1>}(t),\ t>0
$$
holds valid, where $\I_{(\kappa)}^{<n>}$ is the $n$-fold GFI \eqref{GFIn} and $\D_{(k)}^{<n>}$ is the $n$-fold sequential GFD \eqref{GFDLn} in the Riemann-Liouville sense. 

For a function $f\in C_{-1}^{n}(0,+\infty):=\{f:\ f^{(n)}\in C_{-1}(0,+\infty)\}$, the formula
\begin{equation}
\label{sFTLn-C}
(\I_{(\kappa)}^{<n>}\,_*\D_{(k)}^{<n>}\, f) (t) = f(t) - f(0) - \sum_{j=1}^{n-1}\left(\,_*\D_{(k)}^{<j>}\, f\right)(0)\left( \{1\} \, *\, \kappa^{<j>}\right)(t)
\end{equation}
holds valid, where $\I_{(\kappa)}^{<n>}$ is the $n$-fold GFI \eqref{GFIn} and $\,_*\D_{(k)}^{<n>}$ is the $n$-fold sequential GFD \eqref{GFDLn-C}.
\end{theorem}

For the proofs of Theorems \ref{t3-n} and \ref{tgcTfn} and their particular cases we refer the interested readers to \cite{Luc21d}. 

\section{Solutions to the Fractional Differential Equations with the GFDs in the Riemann-Liouville Sense in Terms of the Convolution Series}

First, we introduce the convolution series and treat some of their properties needed for the further discussions. 

For a Sonine kernel $\kappa \in \mathcal{L}_{1}$,  a convolution series in form 
\begin{equation}
\label{conser}
\Sigma_\kappa(t) = \sum_{j=0}^{+\infty} a_j\, \kappa^{<j+1>}(t),\ a_j\in \R
\end{equation}
was introduced in \cite{Luc21c} for analytical treatment of the fractional differential equations with the $n$-fold GFDs of the Caputo type by means of an operational calculus developed for these GFDs. In \cite{Luc21d}, some of the results presented in \cite{Luc21c} were extended to convolution series in form \eqref{conser} generated by any function $\kappa \in C_{-1}(0,+\infty)$ (that is not necessarily a Sonine kernel). 

A very important question regarding convergence of the convolution series \eqref{conser} was answered in \cite{Luc21b,Luc21d}. 

\begin{theorem}[\cite{Luc21d}]
\label{t11}
Let a function $\kappa \in C_{-1}(0,+\infty)$  be represented in the form
\begin{equation}
\label{rep}
\kappa(t) = h_{p}(t)\kappa_1(t),\ t>0,\ p>0,\ \kappa_1\in C[0,+\infty)
\end{equation} 
and the convergence radius of the power  series
\begin{equation}
\label{ser}
\Sigma(z) = \sum^{+\infty }_{j=0}a_{j}\, z^j,\ a_{j}\in \Com,\ z\in \Com
\end{equation}
be non-zero. Then the convolution series 
\eqref{conser}
is convergent for all $t>0$ and defines a function  from the space $C_{-1}(0,+\infty)$. 
Moreover, the series 
\begin{equation}
\label{conser_p}
t^{1-\alpha}\, \Sigma_\kappa(t) = \sum^{+\infty }_{j=0}a_{j}\, t^{1-\alpha}\, \kappa^{<j+1>}(t),\ \ \alpha = \min\{p,\, 1\}
\end{equation}
is uniformly convergent for $t\in [0,\, T]$ for any $T>0$.
\end{theorem}

In what follows, we always assume that the coefficients of the convolution series satisfy the condition that the convergence radius of the corresponding power series is non-zero and thus Theorem \ref{t11} is applicable for these convolution series. 

As an example, let us consider  the geometric series
\begin{equation}
\label{geom}
\Sigma(z) = \sum_{j=0}^{+\infty} \lambda^{j}z^{j},\ \lambda \in \Com,\ z\in \Com.
\end{equation}
For $\lambda \not =0$,  the convergence radius $r$ of this series is equal to $1/|\lambda|$ and thus we can apply  Theorem \ref{t11} that says  that the convolution series generated by a function $\kappa \in C_{-1}(0,+\infty)$ in form
\begin{equation}
\label{l}
l_{\kappa,\lambda}(t) = \sum_{j=0}^{+\infty} \lambda^{j}\kappa^{<j+1>}(t),\ \lambda \in \Com
\end{equation}
is convergent for all $t>0$ and defines a function from the space $C_{-1}(0,+\infty)$.

The convolution series $l_{\kappa,\lambda}$ defined by \eqref{l} plays a very important role in the operational calculus for the GFD of Caputo type developed in \cite{Luc21b}. It provides a far reaching generalization of both the exponential function and the two-parameters Mittag-Leffler function in form \eqref{eqRL_sol}. 

Indeed, let us consider the convolution series \eqref{l} in the case of the kernel function $\kappa=\{ 1\}$. Due to  the formula $\kappa^{<j+1>}(t) = \{ 1\}^{<j+1>}(t) = h_{j+1}(t)$ (see \eqref{2-13}), the convolution series \eqref{l} is reduced to the power series for the exponential function:
\begin{equation}
\label{l-Mic}
l_{\kappa,\lambda}(t) = \sum_{j=0}^{+\infty} \lambda^{j}h_{j+1}(t) =
\sum_{j=0}^{+\infty} \frac{(\lambda\, t)^j}{j!} = e^{\lambda\, t}.
\end{equation}

For the kernel $\kappa(t) = h_{\alpha}(t)$ of the Riemann-Liouville fractional integral,  the formula  $\kappa^{<j+1>}(t) = h_{\alpha}^{<j+1>}(t) = h_{(j+1)\alpha}(t)$ (see \eqref{2-13}) holds valid. Thus, the convolution series \eqref{l} takes the form
\begin{equation}
\label{l-Cap}
l_{\kappa,\lambda}(t) = \sum_{j=0}^{+\infty} \lambda^{j}h_{(j+1)\alpha}(t) =
t^{\alpha-1}\sum_{j=0}^{+\infty} \frac{\lambda^j\, t^{j\alpha}}{\Gamma(j\alpha+\alpha)} = t^{\alpha -1}E_{\alpha,\alpha}(\lambda\, t^{\alpha})
\end{equation}
that is the same as the two-parameters Mittag-Leffler function \eqref{eqRL_sol}. 

For $\kappa \in \mathcal{L}_{1}$, another important convolution series was introduced in  \cite{Luc21b}  as follows:
\begin{equation}
\label{L}
L_{\kappa,\lambda}(t) = (k \, *\ l_{\kappa,\lambda})(t) = 1 + \left(\{ 1 \} * \sum_{j=1}^{+\infty} \lambda^{j}\kappa^{<j>}(\cdot)\right)(t),\ \lambda \in \Com,
\end{equation}
where $k$ is the Sonine kernel associated to the kernel $\kappa$. It is easy to see that in the case $\kappa=\{ 1\}$, the convolution series \eqref{L} coincides with the exponential function: 
\begin{equation}
\label{L-Mic}
L_{\kappa,\lambda}(t) = 1 + \left(\{ 1 \} * \sum_{j=1}^{+\infty} \lambda^{j}h_j(\cdot) \right)(t)=
1+ \sum_{j=1}^{+\infty} \lambda^{j}h_{j+1}(t) = e^{\lambda\, t}.
\end{equation}

In the case of the kernel $\kappa(t) = h_{\alpha}(t),\ t>0,\ 0<\alpha <1$,  the convolution series $L_{\kappa,\lambda}$ is reduced to the two-parameters Mittag-Leffler function \eqref{eqC_sol}:
\begin{equation}
\label{L-Cap-op}
L_{\kappa,\lambda}(t)\, = \, 1 + \left(\{ 1 \} * \sum_{j=1}^{+\infty} \lambda^{j}h_{j\alpha}(\cdot)\right)(t) =
1 + \sum_{j=1}^{+\infty} \lambda^{j}h_{j\alpha+1}(t) = E_{\alpha,1}(\lambda\, t^{\alpha}).
\end{equation}

Analytical solutions to the single- and multi-terms  fractional differential equations with the $n$-fold GFDs of the Caputo type  were presented in \cite{Luc21b} in terms of the convolution series $l_{\kappa,\lambda}$ and $L_{\kappa,\lambda}$. In the rest of this section, we treat the linear single- and multi-terms  fractional differential equations with the $n$-fold GFDs in the Riemann-Liouville sense.

We start with the following auxiliary result:

\begin{theorem}
\label{eqconv}
Two convolution series generated by the same Sonine kernel $\kappa  \in \mathcal{L}_{1}$ coincide for all $t>0$, i.e., 
\begin{equation}
\label{ser12}
\sum_{j=0}^{+\infty} b_{j}\, \kappa^{<j+1>}(t) \equiv  \sum_{j=0}^{+\infty} c_{j}\, \kappa^{<j+1>}(t),\ t>0
\end{equation}
if and only if the corresponding coefficients of these series are equal: 
\begin{equation}
\label{ser13}
a_j = b_j,\ j=0,1,2,\dots .
\end{equation}
\end{theorem}

\begin{proof} 
If the corresponding coefficients of two convolution series generated by the same Sonine kernel $\kappa  \in \mathcal{L}_{1}$ are equal, then we have just one series and evidently the identity \eqref{ser12} holds valid. 

The idea of the proof of the second part of this theorem is the same as the one for the proof of the analogous calculus result for the power series, i.e., under the condition that the identity \eqref{ser12} holds valid we first show that $b_0=c_0$ and then apply the same arguments to prove that $b_1=c_1$, $b_2=c_2$, etc. 

According to Theorem \ref{t11}, the convolution series in form \eqref{conser} is uniformly convergent on any interval $[\epsilon,\ T]$, and thus we can apply the GFI $\I_{(k)}$ to this series term by term:
$$
\left( \I_{(k)}\, \sum_{j=0}^{+\infty} a_{j}\, \kappa^{<j+1>}(\cdot )\right)(t)  = 
 \sum_{j=0}^{+\infty} \left( \I_{(k)}\, a_j\, \kappa^{<j+1>}(\cdot)\right)(t) = 
 $$
 $$
\sum_{j=0}^{+\infty} \left(  a_j\, (k(\cdot)\, *\, \kappa^{<j+1>}(\cdot )\right)(t) = 
a_0 + \sum_{j=1}^{+\infty} a_{j}\, \left(\{ 1\} \, *\, \kappa^{<j>}(\cdot)\right)(t) =  
$$
$$
a_0 + \left( \{1\}\, *\, \sum_{j=1}^{+\infty} a_j\, \kappa^{<j>}(\cdot)\right)(t) = a_0 + (\{1\}\, *\, f_1)(t),
$$
where $f_1$ is the following convolution series:
\begin{equation}
\label{f1}
f_1(t) = \sum_{j=1}^{+\infty} a_j\, \kappa^{<j>}(t)=\sum_{j=0}^{+\infty} a_{j+1}\, \kappa^{<j+1>}(t).
\end{equation}
Summarizing the calculations from above, for the convolution series in form \eqref{conser}, the formula
\begin{equation}
\label{Intcon} 
\left( \I_{(k)}\, \sum_{j=0}^{+\infty} a_{j}\, \kappa^{<j+1>}(\cdot )\right)(t) =  a_0 + \left( \{1\}\, *\, \sum_{j=0}^{+\infty} a_{j+1}\, \kappa^{<j+1>}(\cdot)\right)(t)
\end{equation}
holds valid.

Because the convergence radius of the power series $\Sigma_1(t) = \sum_{j=0}^{+\infty} a_{j+1}\, z^j$ is the same as the convergence radius of the power series $\Sigma(t) = \sum_{j=0}^{+\infty} a_{j}\, z^j$, Theorem  \ref{t11} ensures the inclusion $f_1 \in C_{-1}(0,+\infty)$, where $f_1$ is defined by the formula \eqref{f1}.  As have been shown in \cite{LucGor}, the definite integral of a function from $C_{-1}(0,+\infty)$ is a continuous function on the whole interval $[0,\, +\infty)$ that takes the value zero at the point zero:
\begin{equation}
\label{a0_2} 
\left( \{1\}\, * \, f_1 \right) (x) = (I_{0+}^1\, f_1)(x) \in C[0,\ +\infty),\ \ (I_{0+}^1\, f_1)(0) = 0.
\end{equation}

Now we act with the GFI $\I_{(k)}$ on the equality \eqref{ser12} and apply the formula \eqref{Intcon} to get the relation
\begin{equation}
\label{ser12_1}
b_0 + \left( \{1\}\, *\, \sum_{j=0}^{+\infty} b_{j+1}\, \kappa^{<j+1>}(\cdot)\right)(t) \equiv  c_0 + \left( \{1\}\, *\, \sum_{j=0}^{+\infty} c_{j+1}\, \kappa^{<j+1>}(\cdot)\right)(t),\ t>0.
\end{equation}
Substituting the point $t=0$ into the equality \eqref{ser12_1} and using the formula \eqref{a0_2}, we deduce that $b_0 = c_0$. Now we differentiate the equality \eqref{ser12_1} and get the following identity:
\begin{equation}
\label{ser12_2}
\sum_{j=0}^{+\infty} b_{j+1}\, \kappa^{<j+1>}(t) \equiv  \sum_{j=0}^{+\infty} c_{j+1}\, \kappa^{<j+1>}(t),\ t>0.
\end{equation}
This identity has exactly same structure as the identity \eqref{ser12} from Theorem \ref{eqconv}. Thus we can apply the same arguments as above and derive the relation $b_1 = c_1$. By repeating the same reasoning again and again we arrive at the formula \eqref{ser13} that we wanted to prove.
\end{proof}

Now we are ready to apply the method of convolution series for derivation of solutions to the fractional differential equations with the GFDs and start with the fractional relaxation equation with the GFD of the Riemann-Liouville type:
\begin{equation}
\label{eq-1-1}
( \D_{(k)}\, y)(t)  = \lambda y(t), \ \lambda \in \R,\ t>0. 
\end{equation}

As in the case of the power series, we look for a general solution to this equation in form of a convolution series generated by the Sonine kernel $\kappa$ that is an associated kernel to the kernel $k$ of the GFD from the equation \eqref{eq-1-1}:
\begin{equation}
\label{sol-1-1}
y(t) = \sum_{j=0}^{+\infty} b_j\, \kappa^{<j+1>}(t),\ b_j\in \R.
\end{equation}

To proceed, let as first calculate the image of the convolution series \eqref{sol-1-1} by action of the GFD $\D_{(k)}$: 
$$
( \D_{(k)}\, y)(t)  = \left( \D_{(k)}\, \sum_{j=0}^{+\infty} b_j\, \kappa^{<j+1>}(\cdot)\right)(t) = 
\frac{d}{dt}\left( \I_{(k)}\, \sum_{j=0}^{+\infty} b_j\, \kappa^{<j+1>}(\cdot)\right)(t).
$$

In the proof of Theorem \ref{eqconv} we already calculated the image of the convolution series \eqref{sol-1-1} by action of the GFI $\I_{(k)}$ (formula \eqref{Intcon}). Applying this formula, we arrive at the representation 
\begin{equation}
\label{term1}
( \D_{(k)}\, y)(t)  = \frac{d}{dt}\left( b_0 + \left( \{1\}\, *\, \sum_{j=0}^{+\infty} b_{j+1}\, \kappa^{<j+1>}
(\cdot ) \right)(t)\right) = \sum_{j=0}^{+\infty} b_{j+1}\, \kappa^{<j+1>}(t).
\end{equation}

In the next step, we substitute the right-hand side of \eqref{term1} into the equation \eqref{eq-1-1} and get an equality of two convolution series generated by the same kernel $\kappa$:
$$
\sum_{j=0}^{+\infty} b_{j+1}\, \kappa^{<j+1>}(t) = \sum_{j=0}^{+\infty} \lambda\ b_{j}\, \kappa^{<j+1>}(t),\ t>0.
$$
Application of Theorem \ref{eqconv} to the above identity leads to the following relations for the coefficients of the convolution series \eqref{sol-1-1}:
\begin{equation}
\label{coef1}
b_{j+1} = \lambda\, b_j,\ j=0,1,2,\dots .
\end{equation}
The infinite system \eqref{coef1} of linear equations can be easily solved step by step and we arrive at the explicit solution in form
\begin{equation}
\label{coef2}
b_{j} = b_0\, \lambda^j,\ j=1,2,\dots ,
\end{equation}
where $b_0\in \R$ is an arbitrary constant. Summarizing the arguments presented above, we proved the following theorem:

\begin{theorem}
\label{t-eq1}
The general solution to the fractional relaxation equation  \eqref{eq-1-1} with  the GFD \eqref{FDR-L} in the Riemann-Liouville sense 
can be represented  as follows:
\begin{equation}
\label{eig}
y(t) = \sum_{j=0}^{+\infty} b_0 \, \lambda^j\,  \kappa^{<j+1>}(t) = b_0\, l_{\kappa,\lambda}(t),\ b_0\in \R,
\end{equation}
where $l_{\kappa,\lambda}$ is the convolution series \eqref{l}.
\end{theorem}

\begin{remark}
\label{r1}
The constant $b_0$ in the general solution \eqref{eig} to the equation \eqref{eq-1-1} can be determined from a suitably posed initial condition. The form of this initial condition is prescribed  by Theorem \ref{tgcTfn} (see also the formula \eqref{Intcon}). Indeed, setting $n=1$ in the relation \eqref{sFTLn}, we get the following representation of the projector operator of the GFD \eqref{FDR-L} in the Riemann-Liouville sense:
\begin{equation}
\label{proj1}
(P\, f)(t) = f(t) - (\I_{(\kappa)}\, \D_{(k)}\, f) (t) = \left(\I_{(k)}\, f\right) (0)\kappa(t),\  f\in C_{-1,(k)}^{(1)}(0,+\infty).
\end{equation}
Thus, the initial-value problem 
\begin{equation}
\label{eq-1-1-iv}
\begin{cases}
( \D_{(k)}\, y)(t)  = \lambda y(t), \ \lambda \in \R,\ t>0,\\
\left(\I_{(k)}\, y\right) (0) = b_0
\end{cases}
\end{equation}
has a unique solution given by the formula \eqref{eig}.
\end{remark}

In the case of the Sonine kernel $k(t) = h_{1-\alpha}(t),\ 0<\alpha<1$, the  equation \eqref{eq-1-1} is reduced to the equation \eqref{eqRL} with the Riemann-Liouville fractional derivative and its solution \eqref{eig} is exactly the solution \eqref{eqRL_sol} of the equation \eqref{eqRL} in terms of the two-parameters Mittag-Leffler function (see the formula \eqref{l-Cap}). The initial-value problem \eqref{eq-1-1-iv} takes the well-known form
\begin{equation}
\label{eq-1-1-iv-RL}
\begin{cases}
( D_{0+}^\alpha\, y)(t)  = \lambda y(t), \ \lambda \in \R,\ t>0,\\
\left(I_{0+}^{1-\alpha}\, y\right) (0) = b_0.
\end{cases}
\end{equation}
Its unique solution is given by the formula $y(t) = b_0\, t^{\alpha -1} E_{\alpha, \alpha}(\lambda\, t^\alpha)$. 

Now we proceed with the inhomogeneous equation of type \eqref{eq-1-1}
\begin{equation}
\label{eq-1-2}
( \D_{(k)}\, y)(t)  = \lambda y(t) + f(t), \ \lambda \in \R,\ t>0, 
\end{equation} 
where the source function $f$ is represented in form of a convolution series
\begin{equation}
\label{f}
f(t) = \sum_{j=0}^{+\infty} a_j\, \kappa^{<j+1>}(t),\ a_j\in \R.
\end{equation} 
Again, we look for solutions to the equation \eqref{eq-1-2} in form of the convolution series \eqref{sol-1-1}. Applying exactly the same reasoning as above, we arrive at the following infinite system of linear equations for the coefficients of the convolution series \eqref{sol-1-1}:
\begin{equation}
\label{coef1-1}
b_{j+1} = \lambda\, b_j + a_j,\ j=0,1,2,\dots .
\end{equation}
The explicit form of solutions to this system of equations is as follows:
\begin{equation}
\label{coef2-1}
b_{j} = b_0\, \lambda^j\,  +\,  \sum_{i=0}^{j-1}a_i\, \lambda^{j-i-1}, \ j=1,2,\dots,
\end{equation}
where $b_0\in \R$ is an arbitrary constant. Then the general solution to the equation \eqref{eq-1-2} takes the form
$$
y(t) = b_0\, \kappa(t)+ \sum_{j=1}^{+\infty} \left(b_0\, \lambda^j\,  +\,  \sum_{i=0}^{j-1}a_i\, \lambda^{j-i-1}\right) \, \kappa^{<j+1>}(t) = 
$$
$$
b_0\, \sum_{j=0}^{+\infty} \lambda^j\,\kappa^{<j+1>}(t)  + \sum_{j=1}^{+\infty} \sum_{i=0}^{j-1} a_i\, \lambda^{j-i-1}\, \kappa^{<j+1>}(t).
$$
By direct calculations, we  verify that the second sum in the last formula can be written in a more compact form:
$$
\sum_{j=1}^{+\infty} \sum_{i=0}^{j-1} a_i\,  \lambda^{j-i-1}\, \kappa^{<j+1>}(t) = \sum_{i=0}^{+\infty}  a_i \sum_{j=1}^{+\infty}  \lambda^{j-1}\, \kappa^{<j+i+1>}(t) = (f\, *\, l_{\kappa,\lambda})(t),
$$
where the convolution series $l_{\kappa,\lambda}$ is defined by \eqref{l}. We thus proved the following result:

\begin{theorem}
\label{t-eq2}
The general solution to the inhomogeneous equation \eqref{eq-1-2} has the form
\begin{equation}
\label{eig1}
y(t) = b_0\, l_{\kappa,\lambda}(t) +(f\, *\, l_{\kappa,\lambda})(t),\ b_0\in \R,
\end{equation}
where the convolution series $l_{\kappa,\lambda}$ is defined by \eqref{l}. 

The constant $b_0$ is uniquely determined by the initial condition 
\begin{equation}
\label{ic1}
\left(\I_{(k)}\, y\right) (0) = b_0.
\end{equation}
\end{theorem}

Applying Theorem \ref{t-eq2} to the case of the Riemann-Liouville fractional derivative (kernel $k(t) = h_{1-\alpha}(t),\ 0<\alpha<1$), we obtain the well-known result (\cite{LucSri}):

The unique solution to the initial-value problem 
$$
\begin{cases}
( D_{0+}^\alpha\, y)(t)  = \lambda y(t)+f(t), \ \lambda \in \R,\ t>0,\\
\left(I_{0+}^{1-\alpha}\, y\right) (0) = b_0
\end{cases}
$$
is given by the formula
$$
y(t) = b_0\, t^{\alpha -1} E_{\alpha, \alpha}(\lambda\, t^\alpha) + \int_{0}^t \tau^{\alpha -1} E_{\alpha, \alpha}(\lambda\, \tau^\alpha)\, f(t-\tau)\, d\tau.
$$

Let us now consider a linear inhomogeneous multi-term fractional differential equation with the sequential GFDs \eqref{GFDLn} of the Riemann-Liouville type and with the constant coefficients:
\begin{equation}
\label{eq-1-3}
\sum_{i=0}^n\lambda_i(\D_{(k)}^{<i>}\, y)(t)  =  f(t), \ \lambda_i \in \R,\ i=0,1,\dots,n,\ \lambda_n \not = 0,\  t>0, 
\end{equation} 
where the source function $f$ is represented in form of the  convolution series \eqref{f}. 

As  in the case of the single-term equation \eqref{eq-1-2},  we look for solutions to the multi-term equation \eqref{eq-1-3} in form of the convolution series \eqref{sol-1-1}. First we determine the images of the convolution series \eqref{sol-1-1} by action of the sequential GFDs $\D_{(k)}^{<i>}$, $i=1,2,\dots,n$. For $i=1$, the image is provided by the formula \eqref{term1}. For $i=2,\dots,n$, the formula \eqref{term1} is applied iterative and we arrive at the following result:
\begin{equation}
\label{termi}
( \D_{(k)}^{<i>}\, y)(t)   = \sum_{j=0}^{+\infty} b_{j+i}\, \kappa^{<j+1>}(t),\ i=1,2,\dots,n.
\end{equation}
Now we substitute the convolution series \eqref{sol-1-1}, its images by action of the sequential GFDs $\D_{(k)}^{<i>}$, $i=1,2,\dots,n$ provided by the formula \eqref{termi}, and the convolution series \eqref{f} for the source function into the equation \eqref{eq-1-3} and arrive at the following identity:
$$
 \sum_{i=0}^n\lambda_i\, \left(\sum_{j=0}^{+\infty} b_{j+i}\, \kappa^{<j+1>}(t)\right)  =  \sum_{j=0}^{+\infty} a_{j}\, \kappa^{<j+1>}(t),\ t>0.
$$
Application of Theorem \ref{eqconv} to the above identity leads to the following infinite triangular system of linear equations for the coefficients of the convolution series \eqref{sol-1-1}:
\begin{equation}
\label{coef1-3}
\begin{cases}
\lambda_0\, b_0 +\lambda_1\, b_1 +\dots + \ \lambda_n\, b_n = a_0,\\
\lambda_0\, b_1 +\lambda_1\, b_2 +\dots + \ \lambda_n\, b_{n+1} = a_1,\\
\dots \\
\lambda_0\, b_n +\lambda_1\, b_{n+1} +\dots + \ \lambda_n\, b_{2n} = a_n,\\
\lambda_0\, b_{n+1} +\lambda_1\, b_{n+2} +\dots + \ \lambda_n\, b_{2n+1} = a_{n+1} \\
\dots
\end{cases}
\end{equation}
In this system, the first $n$ coefficients ($b_0,\ b_1,\dots,b_{n-1}$) can be chosen arbitrary and all other coefficients are determined step by step as solutions to the infinite triangular system \eqref{coef1-3} of linear equations:
\begin{equation}
\label{coef1-4}
b_{n+l}=(a_l - \lambda_0\, b_l- \dots - \lambda_{n-1} b_{n+l-1})/\lambda_n,\ l=0,1,2,\dots
\end{equation}

We thus proved the following result:

\begin{theorem}
\label{t-eq3}
The general solution to the inhomogeneous multi-term fractional differential equation \eqref{eq-1-3} can be represented as the convolution series \eqref{sol-1-1}, where the first $n$ coefficients ($b_0,\ b_1,\dots,b_{n-1}$) are arbitrary real constants and other coefficients are calculated according to the formula \eqref{coef1-4}.
\end{theorem}

The constants $b_0,\ b_1,\dots,b_{n-1}$ in the general solution to the equation \eqref{eq-1-3} presented in Theorem \eqref{t-eq3}  can be determined based on the suitably posed initial conditions.  The form of these initial conditions is prescribed by Theorem \ref{tgcTfn}. Indeed,  for a function $f\in C_{-1,(k)}^{(n)}(0,+\infty)$, the formula \eqref{sFTLn} can be rewritten as follows: 
\begin{equation}
\label{projn}
(P\, f)(t) = f(t) -  (\I_{(\kappa)}^{<n>}\, \D_{(k)}^{<n>}\, f) (t) = \sum_{j=0}^{n-1}\left( \I_{(k)}\, \D_{(k)}^{<j>}\, f\right)(0)\kappa^{<j+1>}(t),\  t>0,
\end{equation}
where $P$ is the projector operator of the $n$-fold sequential GFD of the Riemann-Liouville type. Thus, to uniquely determine the constants $b_0,\ b_1,\dots,b_{n-1}$ in the general solution, the equation \eqref{eq-1-3} has to be equipped with the initial conditions in the form
\begin{equation}
\label{icm}
\left( \I_{(k)}\, \D_{(k)}^{<j>}\, y\right)(0) = b_j,\ j=0,1,\dots,n-1.
\end{equation}

Finally, we mention that the inhomogeneous multi-term fractional differential equation of type \eqref{eq-1-3} with the sequential Riemann-Liouville fractional derivatives (the case of the kernel $k(t)=h_{1-\alpha}(t)$ in the equation \eqref{eq-1-3}) was treated in \cite{LucSri,Luc99} my using operational calculus of the Mikusi\'nski type for the Riemann-Liouville fractional derivative.





\end{document}